%soli.tex:  In how many ways can you play Stanley Solitaire?
%%a Plain TeX file  Shalosh B. Ekhad and Doron Zeilberger ( pages)

%begin macros

\baselineskip=14pt
\parskip=10pt
\def\halmos{\hbox{\vrule height0.15cm width0.01cm\vbox{\hrule height
  0.01cm width0.2cm \vskip0.15cm \hrule height 0.01cm width0.2cm}\vrule
  height0.15cm width 0.01cm}}

\magnification=\magstephalf

\def\1{{\overline{1}}}
\def\2{{\overline{2}}}
\parindent=0pt
\overfullrule=0in

\def\frac#1#2{{#1 \over #2}}
%\headline={\rm  \ifodd\pageno  \RightHead  \else  \LeftHead  \fi}
%\def\RightHead{\centerline{
%Title
%}}
%\def\LeftHead{ \centerline{Doron Zeilberger}}
%end macros
\centerline
{\bf 
In how many ways can you play Stanley Solitaire?
}
\bigskip
\centerline
{\it Shalosh B. EKHAD and Doron ZEILBERGER}
\bigskip
\qquad {\it Dedicated to Richard Peter Stanley (b. June 23, 1944), the non-bossy superboss of Algebraic Combinatorics, on his $(C_{5}+2 C_{4}+2 C_{3})^{th}$ birthday}

{\bf Preface: The Hottest Problem in Algebraic Combinatorics in the Early 1980s}

When Richard Stanley was half his age (today), the hottest problem in algebraic combinatorics was to
find, and prove, a formula for the number of ways of writing the longest permutation,
$$
n,n-1, \dots , 1 \quad
$$
as a product of $n(n-1)/2$ {\it adjacent transpositions}, $(i,i+1)$ ($1\leq i \leq n-1)$ that famously
generate the symmetric group. Using data communicated by Paul Edelman, Stanley conjectured that these numbers
are the same as the number of so-called Young tableaux of shape $[n-1,n-2, \dots, 1]$.
This was brilliantly proved (in much greater generality) in [St]. Independently it was also proved by Paul Edelman and Curtis Greene [EG],
and outlined by Alain Lascoux and Marcel-Paul Sch\"utzenberger [LSc]. 

Here we will see how {\it highbrow} meets {\it lowbrow}, in a very nice way.

{\bf Stanley Solitaire: A Game You Can Teach a Five-Year-Old}

It is not wise to give candy to children unless they learn something from it. Here is a simple educational game that teaches children
the notions of {\it larger} and {\it smaller}, and {\it left} and {\it right}. We named this game {\it Stanley Solitaire} for a reason to be made clear soon.

You start with an {\it initial position} on a one-dimensional {\bf board}, where some squares have piles of candies (for example M\&Ms) on them.
The left side of the board is fixed, but one can move as far right as one wishes, following the rules.

One such starting position (out of infinitely many) could be
$$
[4,5,0,0,2,0,3,1] \quad,
$$
meaning that the leftmost location has $4$ candies, the second-from-the-left has $5$ candies, the third, fourth, and sixth locations are empty, the fifth has $2$ candies, the seventh has $3$ candies, and the eighth has
one candy.

The player is allowed to eat all the candies, provided that he or she follow the following rule of what constitutes  {\bf one legal move}.

{\it Choose any two {\bf adjacent} piles where the number of candies in the left pile is (strictly) larger than the the number of candies in the (possibly empty) right pile, exchange them, and eat one of the candies of
the larger pile (the one that initially was on the left, and now is to the right)}. 

The game continues until all the candies are eaten.

As soon as the player messes up, they are not allowed to finish the game, so they have strong motivation to follow the rules.

In the above example, the legal moves are

$\bullet$ Eat one candy from the second location, and exchange the second and third,  getting to the new position $[4,0,4,0,2,0,3,1]$ \quad .

$\bullet$ Eat one candy from the fifth location, and exchange the fifth and the sixth piles, getting to the new position $[4,5,0,0,0,1,3,1]$ \quad .

$\bullet$ Eat one candy from the seventh location, and exchange the seventh  and eighth piles, getting to the new position $[4,5,0,0,2,0,1,2]$ \quad .

$\bullet$ Eat one candy from the eighth location, and exchange the eighth and the (implied)  ninth pile getting to the new position $[4,5,0,0,2,0,3,0,0]$, but
we delete all rightmost zeroes, so the new position is $[4,5,0,0,2,0,3]$.

More formally, for a position
$$
[a_1, \dots, a_k] \quad,
$$
where we insist that $a_1>0$ and  $a_k>0$, and otherwise $a_i \geq 0$, the legal moves are as follows.

For $1 \leq i <k$, if $a_i>a_{i+1}$ then go to position

$$
[a_1, \dots, a_{i-1} \,, \, a_{i+1}, a_i -1 \,, \, a_{i+2} \,,  \,\dots, a_k] \quad,
$$
and yet another legal move is to position $[a_1, \dots, a_{k-1},0,  a_k-1]$. Whenever the new position has $0$ either on the left or right, we remove them, keeping the convention that the leftmost and rightmost piles are non-empty.

{\bf A problem you can Explain To a Six-Year-Old}

{\it For any starting position, in how many ways can you legally eat all the candies?}

For example, if the initial position is $[2,1]$, then there are only two ways:

$$
[2,1] \rightarrow [1,1] \rightarrow [1] \rightarrow [] \quad ,
$$
$$
[2,1] \rightarrow [2,0,0] (\,alias \, [2]\,) \rightarrow [1] \rightarrow [] \quad ,
$$

while if the starting position is $[2,2,1]$ there are five ways. Here they are:
$$
[2, 2, 1], [2, 1, 1], [1, 1, 1], [1, 1], [1], [] \quad,
$$
$$
[2, 2, 1], [2, 1, 1], [2, 1], [1, 1], [1], [] \quad,
$$
$$
[2, 2, 1], [2, 1, 1], [2, 1], [2], [1], [] \quad,
$$
$$
[2, 2, 1], [2, 2], [2, 0, 1], [1, 1], [1], [] \quad,
$$
$$
[2, 2, 1], [2, 2], [2, 0, 1], [2], [1], [] \quad.
$$

{\bf An Answer you can Explain to a Seven-Year-Old}

Thanks to the seminal work of Richard Stanley [S], and independently, Paul Edelman, and Curtis Greene [EG] (and presumably [LSc], but it is over our heads), there is
an extremely easy and elegant {\bf answer}, that only involves addition and multiplication (recall that $a!=1 \cdot 2 \cdots a$).

{\bf Theorem 1}: The {\bf exact} number of ways to play Stanley-Solitaire starting with position
$$
[a_1, \dots, a_k] \quad,
$$
where $a_1 \geq  a_2 \geq \dots \geq a_k >0$  equals:
$$
\frac{(a_1+ \dots + a_k)!}{(a_1+k-1)! (a_2+k-2)! \cdots (a_k)!} \cdot \prod_{1 \leq i < j \leq k} (a_i-a_j+j-i) \quad .
\eqno(YFM)
$$

This simply stated theorem, that {\it begs} for a simple proof, is surprisingly deep. We will use:

{\bf Stanley's Amazing Theorem} ([St], Cor. 4.2,  p. 368): If $w$ is a  permutation for which the two partitions  
$\lambda(w)$ and $\mu(w)$ (see [St], p. 367, lines 3-7 for their definitions)  are both equal to $[a_1, \dots, a_k]$, 
then the number of ways of playing Stanley Solitaire on  that initial configuration is given  by $(YFM)$.

Let's use this deep theorem to prove our own Theorem 1.

{\bf Proof of  that Stanley's Theorem Implies Theorem 1}: 
If $a_1>a_2 >\dots > a_k \geq 1$, define the permutation of $\{1,2, \dots, a_1+1\}$ as follows:
$$
w:=[a_1+1, a_2+1, \dots, a_k+1, 1,2, \dots, a_k, a_k+2, \dots, a_{k-1},a_{k-1}+2, \dots, a_2, a_{2}+2, \dots , a_1] \quad,
$$
that has the property that $\lambda(w)=\mu(w)=[a_1, \dots, a_k]$ (check!) \halmos .

{\bf Why is Our Theorem 1 So Hard to Prove Directly?}

Even the case of two initial piles of candy, $[a_1,a_2]$, with $a_1>a_2$ is surprisingly hard to prove. The natural way would be to use (symbolic) dynamical programming.
Let $a_1 \geq a_2 >0$, and let $S([a_1,a_2])$ be the number of ways of playing starting with $[a_1,a_2]$.

We have:
$$
S([a_1,a_2]) \, = \, S([a_2,a_1-1]) + S([a_1,0,a_2-1]) \quad,
$$
so we are {\bf forced} to try and find explicit expressions for two new quantities, $S([a_2,a_1])$ with $a_2 < a_1$, and $S([a_1,0,a_2])$, with $a_1>a_2$.
Regarding the first one we have
$$
S([a_2,a_1])= S([a_2,0,a_1-1] \quad,
$$
and the second one,
$$
S([a_1,0,a_2])= S([a_1-1,a_2])+ S([a_1,0,0,a_2-1]) \quad.
$$
These new quantities are no longer `nice', and while it is possible to conjecture closed-form expressions for each specific template, of the form $S([a_1, 0^{k}, a_2])$, for $a_1>a_2$ and
$S([a_2, 0^{k}, a_1])$, they get nastier and nastier as $k$ gets larger, and we go down a {\it slippery slope}.

The conjectured expressions can be gotten either by {\it guessing} (see Maple package below), or
by {\it going backwards}, trusting the theorem, and using the implied recurrences backwards, and recursively.
But this would be {\it cheating}, and eventually we would need to check the {\it base cases}. We were unable to carry this plan,
and this made us admire [St] even more. So it takes a very smart $24$-year graduate student who read {\it Enumerative Combinatorics}, vol. 2,
to {\bf understand} {\it why}  is the Stanley Solitaire theorem true.

{\bf Comment}: More generally, for any $231$-avoiding permutation $\pi=\pi_1 \dots \pi_k$ (there are $C_k=(2k)!/(k!(k+1)!)$ of them) and $a_1 \geq a_2 \geq a_k >0$,
the number of ways of playing Stanley Solitaire with initial position $[a_{\pi_1}, \dots, a_{\pi_k}]$ is also given by $(YFM)$. 

If $\pi$ is not $231$-avoiding, then the formulas are no longer so nice (see below for an example)).

{\bf The Maple package} {\tt StanleySolitaire.txt}

You are welcome to download, and experiment with, the Maple package

{\tt http://www.math.rutgers.edu/\~{}zeilberg/tokhniot/StanleySolitaire.txt}  \quad,

that has the procedure {\tt S(a)}, that uses {\it numerical dynamical programming} to compute these quantities. This would enable the seven-year-old to {\it conjecture} Theorem 1, for at least $S([a_1,a_2])$, with $a_1 \geq a_2  >0$.

It can also guess (increasingly complicated) closed-form expressions for other scenarios. For example for $S([b,c,a])$, for which $(YFM)$ is no longer valid, we have instead the following formula.

{\bf Fact}: The number of ways of playing Stanley Solitaire starting with initial position $[b,c,a]$, where $a \geq b \geq c >0$ is:
$$
\frac{(a+b+c)!}{(a+3)!(b+1)! c!} \cdot (a^{2} b -a^{2} c +a \,b^{2}-2 a b c +a \,c^{2}-b^{2} c +b \,c^{2}+a^{2}+5 a b -6 a c +3 b^{2}-6 b c +3 c^{2}+4 a +6 b -9 c +3) (a -b +2) \quad .
$$

Since this formula is obviously true, we didn't bother to prove it, but readers are most welcome to attempt a proof, either directly (that seems very hard) or by using [St], like we did for Theorem 1.

The front of this article  

{\tt https://sites.math.rutgers.edu/\~{}zeilberg/mamarim/mamarimhtml/soli.html}  \quad ,

contains numerous sample input and output files, but you are welcome to use our package to generate more data, and who knows?, find an elementary, {\it high school algebra}, proof of the theorem, that
a nine-year-old can understand (and enjoy),  without the heavy machinery of algebraic combinatorics.

{\bf Conclusion}: We fully enjoyed attending Stanley@60,  Stanley@70,  and Stanley@80. We are  looking forward to Stanley@90.

{\bf Acknowledgment}: Many thanks to ex-president Phil Hanlon, yet-another non-bossy {\it superboss}, who inspired the dedication.

{\bf References}

[EG] Paul H. Edelman and Curtis Greene, {\it Balanced tableaux}, Adv.  Math. {\bf 63} (1987), 42-99. \hfill\break
{\tt https://core.ac.uk/download/pdf/82733311.pdf}

[LSc] Alain Lascoux and Marcel-Paul Sch\"utzenberger, {\it Structure de Hopf de l'anneau de cohomologie et de l'anneau de Grothendick d'une
vari\'et\'e de drapeaux}, C.R. Acad. Sc. Paris {\bf 295}, S\'erie 1, 629-633. \hfill\break
{\tt http://www-igm.univ-mlv.fr/\~{}berstel/Mps/Travaux/A/1982-2HopfCras.pdf}

[St] Richard P. Stanley, {\it On the number of reduced decompositions of elements of the Coxeter group},
Europ. J. Combinatorics {\bf 5} (1984), 359-372. \hfill\break
{\tt  https://math.mit.edu/\~{}rstan/pubs/pubfiles/56.pdf}

\bigskip
\hrule
\bigskip

Shalosh B. Ekhad and Doron Zeilberger, Department of Mathematics, Rutgers University (New Brunswick), Hill Center-Busch Campus, 110 Frelinghuysen
Rd., Piscataway, NJ 08854-8019, USA. \hfill\break
Email: {\tt ShaloshBEkhad at gmail  dot com}, {\tt DoronZeil at gmail  dot com}   \quad .
\bigskip
June 23, 2024.

\end